\documentclass[11pt,a4paper]{article}
\usepackage{a4,amssymb,latexsym,amsmath,color}
\begin{document}


\newcounter{z}
\newenvironment{enum}{\setcounter{z}{0}
\begin{list}{\rm\alph{z})}{\usecounter{z}
\setlength{\topsep}{1ex}\setlength{\labelwidth}{0.5cm}
\setlength{\leftmargin}{1cm}\setlength{\labelsep}{0.25cm}
\setlength{\parsep}{-3pt}}}{\end{list}~\\[-6ex]}

\newcounter{z1}
\newenvironment{enum1}{\setcounter{z1}{0}
\begin{list}{\arabic{z1})}{\usecounter{z1}
\setlength{\topsep}{1ex}\setlength{\labelwidth}{0.6cm}
\setlength{\leftmargin}{1cm}\setlength{\labelsep}{0.25cm}
\setlength{\parsep}{-3pt}}}{\end{list}~\\[-6ex]}

\newcounter{z2}
\newenvironment{enum2}{\setcounter{z2}{0}
\begin{list}{$\bullet$}{\usecounter{z2}
\setlength{\topsep}{1ex}\setlength{\labelwidth}{0.6cm}
\setlength{\leftmargin}{0.75cm}\setlength{\labelsep}{0.25cm}
\setlength{\parsep}{-3pt}}}{\end{list}~\\[-6ex]}

\newcounter{zb}
\newenvironment{enumbib}{\setcounter{zb}{0}
\begin{small}
\begin{list}{\arabic{zb}.}{\usecounter{zb}
\setlength{\topsep}{1ex}\setlength{\labelwidth}{0.6cm}
\setlength{\leftmargin}{0.8cm}\setlength{\labelsep}{0.25cm}
\setlength{\parsep}{1pt}}}{\end{list}\end{small}}

\newcounter{pic}\setcounter{pic}{0}

\newcommand{\bdpm}{\begin{displaymath}}
\newcommand{\edpm}{\end{displaymath}}

\newcommand{\beas}{\begin{eqnarray*}}
\newcommand{\eeas}{\end{eqnarray*}}

\newcommand{\ba}{\begin{array}}
\newcommand{\ea}{\end{array}}

\newenvironment{proof}{\vspace*{-0.25cm}\begin{sloppypar}\noindent{\it 
Proof.}}{\hfill\qed\end{sloppypar}\vspace*{0.15cm}}

\newtheorem{theo}{Theorem}[section]          
\newtheorem{lem}[theo]{Lemma}
\newtheorem{cor}[theo]{Corollary}
\newtheorem{prop}[theo]{Proposition}
\newtheorem{exa}[theo]{Example}
\newtheorem{rem}[theo]{Remark}
\newtheorem{rems}[theo]{Remarks}

\newcommand{\brm}{\begin{rm}}
\newcommand{\erm}{\end{rm}}

\newcommand{\qed}{\hfill $\Box$}
\newcommand{\lis}{{\sf lis}}
\newcommand{\sg}{{\sf sign}}
\newcommand{\srs}{{\sf srs}}
\newcommand{\inv}{{\sf inv}}
\newcommand{\ldes}{{\sf ldes}}
\newcommand{\lind}{{\sf lind}}


\begin{center}
{\large\bf REFINED SIGN-BALANCE ON 321-AVOIDING PERMUTATIONS}\\[1cm]
Astrid Reifegerste\\
Institut f\"ur Mathematik, Universit\"at Hannover\\
Welfengarten 1, 30167 Hannover, Germany\\
{\it reifegerste@math.uni-hannover.de}\\[0.5cm]
version of May 22, 2003 
\end{center}
\vspace*{0.5cm}

\begin{footnotesize}
{\sc Abstract.} The number of even $321$-avoiding permutations of length $n$ is equal to the number 
of odd ones if $n$ is even, and exceeds it by the $\frac{n-1}{2}$th Catalan number 
otherwise. We present an involution that proves a refinement of this 
sign-balance property respecting the length of the longest increasing 
subsequence of the permutation. In addition, this yields a combinatorial proof 
of a recent analogous result of Adin and Roichman dealing with the last 
descent. In particular, we answer the question how to obtain the sign of a 
$321$-avoiding permutation from the pair of tableaux resulting from the 
Robinson-Schensted-Knuth algorithm. The proof of the simple solution bases on a 
matching method given by Elizalde and Pak.
\end{footnotesize}
\vspace*{1cm}


\setcounter{section}{1}\setcounter{theo}{0}

\centerline{\large{\bf 1}\hspace*{0.25cm}
{\sc Introduction}}
\vspace*{0.5cm}

Let $T_n$ be the set of $321$-avoiding permutations in the symmetric group 
${\cal S}_n$. (A permutation is called {\it $321$-avoiding} if it has no 
decreasing subsequence of length three.) Simion and Schmidt 
\cite{simion-schmidt} proved the following sign-balance property of $T_n$: the number of even permutations in $T_n$ is equal to the number 
of odd permutations if $n$ is even, and exceeds it by the Catalan number 
$C_{\frac{1}{2}(n-1)}$ otherwise. Very recently, Adin and Roichman \cite{adin-roichman} refined this result 
by taking into account the maximum descent. We give an analogous result for a further important 
permutation statistic, the length of the longest increasing subsequence. In a 
recent paper, Stanley \cite{stanley1} established the importance of the 
sign-balance.\\[2ex] 
For $\pi\in T_n$, let $\lis(\pi)$ be the length of the longest increasing 
subsequence in $\pi$. By \cite[Thm. 4]{deutsch etal}, the number of permutations $\pi\in T_n$ for which $\lis(\pi)=k$ is just the square 
of the {\it ballot number}
\bdpm
b(n,k)=\frac{2k-n+1}{n+1}{n+1\choose k+1}
\edpm
where $\lfloor\frac{n+1}{2}\rfloor\le k\le n$. This is an immediate consequence of the Robinson-Schensted-Knuth correspondence 
which gives a bijection between permutations and pairs of standard Young tableaux of the same shape 
(see, e.g., \cite{stanley}). It is well-known that the 
length of the longest increasing (decreasing) subsequence of a permutation is 
just the length of the first row (column) of its associated tableaux.\\[2ex]
Consequently, $321$-avoiding permutations correspond to pairs of 
standard Young tableaux having at most two rows. Such tableaux can be identified 
with ballot sequences in a natural way. (We call a sequence $b_1b_2\cdots b_n$ with $b_i=\pm1$ 
a {\it ballot sequence} if $b_1+b_2+\ldots+b_j\ge0$ for all $j$.) 
Given a standard Young tableau with at most two rows, set $b_i=1$ if $i$ appears in the first 
row and $b_i=-1$ otherwise. By the Ballot theorem (see, e.g., \cite{feller}), there 
are $b(n,k)$ ballot sequences of length $n$ having exactly $k$ components 
equal $1$.\\
Elizalde and Pak \cite{elizalde-pak} already made use of the simplicity of the 
RSK correspondence for $321$-avoiding permutations for their bijection between 
refined restricted permutations.\\[2ex]  
Our main result is the following. 
\begin{theo} \label{main theorem}
For all $n\ge1$, we have
\beas
\sum_{\pi\in T_{2n+1}} \sg(\pi)\cdot q^{\lis(\pi)}&=& 
\sum_{\pi\in T_n} q^{2\lis(\pi)+1}\\[1ex]
\sum_{\pi\in T_{2n+2}} \sg(\pi)\cdot q^{\lis(\pi)}&=&
(q-1)\sum_{\pi\in T_n} q^{2\lis(\pi)+1}.
\eeas
\end{theo}
To prove this theorem, we establish an involution on $T_n$ which is defined in 
terms of the corresponding tableau pairs. While it is clear how to see the length of 
the longest increasing subsequence from the tableaux, for the sign of a 
permutation this connection is unknown until now. In Section 2, we solve this problem for $321$-avoiding permutations.\\ 
The description of the main bijection $\Phi$ is done in Section 3. Modifying $\Phi$ slightly yields a combinatorial proof of Adin-Roichman's result 
which is given in Section 4. We conclude with a simple combinatorial proof of 
the equidistribution of the last descent and $\pi^{-1}(n)-1$ over $T_n$, Adin and Roichman asked for.
\vspace*{1cm}


\setcounter{section}{2}\setcounter{theo}{0}

\centerline{\large{\bf 2}\hspace*{0.25cm}
{\sc How to obtain the sign from the tableaux}}
\vspace*{0.5cm}

The key problem we are confronted with is to figure out how to see 
the sign of a permutation by looking at its associated pair of tableaux. 
For 321-avoiding permutations, there is a simple answer.\\[2ex] 
Given $\pi\in T_n$, let $p$ and $q$ be the ballot sequences defined by the 
tableaux $P$ and $Q$ which we obtain by applying the RSK algorithm to $\pi$. 
(As usual, we write $P$ to denote the insertion tableau, and $Q$ for the 
recording tableau.)
Define the statistic $\srs$ to be the sum of the elements of the second 
row of $P$ and $Q$, that is,
\bdpm
\srs(\pi)=\sum_{p_i=-1}i+\sum_{q_i=-1}i.
\edpm

\begin{prop} \label{sign}
For any $\pi\in T_n$ with $\lis(\pi)=k$, we have $\sg(\pi)=(-1)^{\srs(\pi)+n-k}$.
\end{prop}

For the proof we use a method to generate a matching between excedances and 
anti-excedances described by Elizalde and Pak \cite{elizalde-pak}.\\ 
For $\pi\in{\cal S}_n$, an {\it excedance} ({\it anti-excedance}) of $\pi$ is an integer $i$ 
for which $\pi_i>i$ ($\pi_i<i$). Here the element $\pi_i$ is called an {\it 
excedance letter} ({\it anti-excedance letter}). It is characteristic for $321$-avoiding 
permutations that both the subword consisting of all excedance letters and the 
subword consisting of the remaining letters are increasing. Due to this 
condition, we say that the permutation is {\it bi-increasing}. 
In particular, the fixed points of a permutation $\pi\in T_n$ are just such integers $i$ for which 
$\pi_j<\pi_i$ for all $j<i$ and $\pi_j>\pi_i$ for all $j>i$. Consequently, each 
longest increasing subsequence of $\pi$ contains all the fixed points of 
$\pi$.\\[2ex]
Given a permutation $\pi\in T_n$, let $i_1<i_2<\ldots<i_s$ be its excedances and $j_1<j_2<\ldots<j_t$ the 
anti-excedances. Construct the matching as follows.  
\begin{enum1} 
\item Initialize $a=b=1$.
\item While $a\le s$ and $b\le t$ repeat the procedure:
\begin{enum2}
\item If $i_a>j_b$, then increase $b$ by $1$.
\item If $\pi_{i_a}<\pi_{j_b}$, then increase $a$ by $1$.
\item If $i_a<j_b$ and $\pi_{i_a}>\pi_{j_b}$, then match $i_a$ with $j_b$ and 
increase both $a$ and $b$ by $1$.
\end{enum2}
\end{enum1}
\vspace*{-0.35cm}

For more clarity, we represent $\pi\in T_n$ by an $n\times n$ array with a dot in each of the squares $(i,\pi_i)$. We will identify the 
integer $i$ with the dot $(i,\pi_i)$, and use the terms ``excedance" and 
``anti-excedance" correspondingly for the dots as well.
\vspace*{2ex}

\begin{center}
\addtocounter{pic}{1}
\unitlength=0.35cm
\begin{picture}(12,12)
\linethickness{0.3pt}
\multiput(0,0)(0,1){13}{\line(1,0){12}}
\multiput(0,0)(1,0){13}{\line(0,1){12}}
\bezier{100}(0,12)(6,6)(12,0)
\put(3.5,11.5){\circle*{0.4}}
\put(0.5,10.5){\circle*{0.4}}
\put(1.5,9.5){\circle*{0.4}}
\put(4.5,8.5){\circle*{0.4}}
\put(6.5,7.5){\circle*{0.4}}
\put(7.5,6.5){\circle*{0.4}}
\put(2.5,5.5){\circle*{0.4}}
\put(5.5,4.5){\circle*{0.4}}
\put(8.5,3.5){\circle*{0.4}}
\put(11.5,2.5){\circle*{0.4}}
\put(9.5,1.5){\circle*{0.4}}
\put(10.5,0.5){\circle*{0.4}}
\linethickness{1pt}
\bezier{100}(3.5,11.5)(2,11)(0.5,10.5)
\bezier{100}(4.5,8.5)(3.5,7)(2.5,5.5)
\bezier{100}(6.5,7.5)(6,6)(5.5,4.5)
\bezier{100}(11.5,2.5)(10.5,2)(9.5,1.5)
\end{picture}

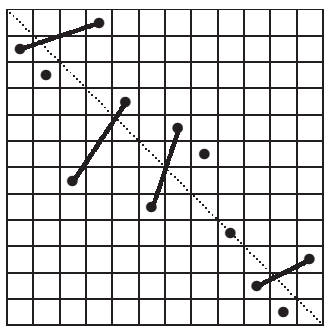
{\footnotesize{\bf Figure \thepic}\hspace*{0.25cm}Example of the matching for 
$\pi=4\:1\:2\:5\:7\:8\:3\:6\:9\:12\:10\:11\in T_{12}$}
\end{center}

It follows immediately from the descriptions of the RSK algorithm and the matching that the element $\pi_i$ appears in the second row of the tableau $P$ if and only 
if $(i,\pi_i)$ is a matched excedance. In particular, the number of matched 
pairs is just the length of the second row of $P$, see \cite[Lemma 
5]{elizalde-pak}. Because of the symmetry of the RSK correspondence -- if 
$(P,Q)$ is the pair of tableaux associated with $\pi\in{\cal S}_n$, then $(Q,P)$ corresponds 
to the inverse permutation $\pi^{-1}$ -- the integer $j$ is contained in the 
second row of $Q$ if and only if $(j,\pi_j)$ is a matched anti-excedance.\\
For the permutation $\pi\in T_{12}$ considered in the above figure, we obtain the 
tableaux 
\vspace*{1ex}
\bdpm
P=
\mbox{\begin{tabular}{|p{1mm}|p{1mm}|p{1mm}|p{1mm}|p{1mm}|p{1mm}|p{1mm}|p{1mm}|p{1mm}|p{1mm}|p{1mm}|p{1mm}|}
\cline{1-8}
\hspace*{-0.2mm}{\footnotesize1}&\hspace*{-0.2mm}{\footnotesize2}&
\hspace*{-0.2mm}{\footnotesize3}&\hspace*{-0.2mm}{\footnotesize6}&
\hspace*{-0.2mm}{\footnotesize8}&\hspace*{-0.2mm}{\footnotesize9}&
\hspace*{-1mm}{\footnotesize10}&\hspace*{-1mm}{\footnotesize11}\\\cline{1-8}
\hspace*{-0.2mm}{\footnotesize4}&\hspace*{-0.2mm}{\footnotesize5}&
\hspace*{-0.2mm}{\footnotesize7}&\hspace*{-1mm}{\footnotesize12}&\multicolumn{8}{|c}{}\\\cline{1-4}
\end{tabular}
\hspace*{3cm}}
Q=\;
\mbox{\begin{tabular}{|p{1mm}|p{1mm}|p{1mm}|p{1mm}|p{1mm}|p{1mm}|p{1mm}|p{1mm}|p{1mm}|p{1mm}|p{1mm}|p{1mm}|}
\cline{1-8}
\hspace*{-0.2mm}{\footnotesize1}&\hspace*{-0.2mm}{\footnotesize3}&
\hspace*{-0.2mm}{\footnotesize4}&\hspace*{-0.2mm}{\footnotesize5}&
\hspace*{-0.2mm}{\footnotesize6}&\hspace*{-0.2mm}{\footnotesize9}&
\hspace*{-1mm}{\footnotesize10}&\hspace*{-1mm}{\footnotesize12}\\\cline{1-8}
\hspace*{-0.2mm}{\footnotesize2}&\hspace*{-0.2mm}{\footnotesize7}&
\hspace*{-0.2mm}{\footnotesize8}&\hspace*{-1mm}{\footnotesize11}&\multicolumn{8}{|c}{}\\\cline{1-4}
\end{tabular}\hspace*{2cm}}
\edpm
\vspace*{-0.2cm}

and hence $\srs(\pi)=56$.\\[2ex]
The connection between the matched pairs and the inversions of the 
permutation is not difficult to see. We call a pair $(i,j)$ an {\it inversion} of $\pi$ if $i<j$ and 
$\pi_i>\pi_j$, and write $\inv(\pi)$ to denote the number of inversions of 
$\pi$.\\
Let $[i,j]$ be a matched pair where $i<j$. 
Denote by $c(i,j)$ the number of integers $l$ for which $i<l<j$ and $\pi_l>\pi_i$ 
or $j<l$ and $\pi_j<\pi_l<\pi_i$. Graphically, the dots $(l,\pi_l)$ are  
contained in the regions 2 and 3 of the scheme: 
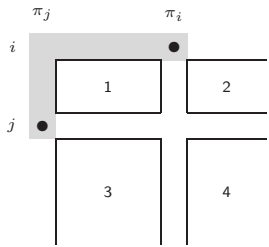
\begin{center}
\unitlength=0.35cm
\definecolor{gray1}{gray}{0.85}
\fboxsep0cm
\fboxrule0cm
\begin{picture}(9,8.5)
\put(0,7){\fcolorbox{gray1}{gray1}{\makebox(6,1){}}}
\put(0,4){\fcolorbox{gray1}{gray1}{\makebox(1,4){}}}
\put(5.5,7.5){\circle*{0.4}}
\put(0.5,4.5){\circle*{0.4}}
\linethickness{0.3pt}
\multiput(1,5)(0,2){2}{\line(1,0){4}}
\multiput(1,5)(4,0){2}{\line(0,1){2}}
\put(1,4){\line(1,0){4}}
\multiput(1,0)(4,0){2}{\line(0,1){4}}
\multiput(6,5)(0,2){2}{\line(1,0){3}}
\put(6,5){\line(0,1){2}}
\put(6,4){\line(1,0){3}}
\put(6,0){\line(0,1){4}}
\put(3,6){\makebox(0,0)[cc]{\sf\tiny1}}
\put(7.5,6){\makebox(0,0)[cc]{\sf\tiny2}}
\put(3,2){\makebox(0,0)[cc]{\sf\tiny3}}
\put(7.5,2){\makebox(0,0)[cc]{\sf\tiny4}}
\put(-0.5,7.5){\makebox(0,0)[cr]{\tiny$i$}}
\put(-0.5,4.5){\makebox(0,0)[cr]{\tiny$j$}}
\put(5.5,8.5){\makebox(0,0)[cb]{\tiny$\pi_i$}}
\put(0.5,8.5){\makebox(0,0)[cb]{\tiny$\pi_j$}}
\end{picture}
\end{center}
Clearly, the first mentioned ones are excedances while the last mentioned are 
anti-excedances. For $r=1,\ldots,4$, let $c_r(i,j)$ be the number of dots in the region $r$.

\begin{lem}
Let $\pi\in T_n$ with $\lis(\pi)=k$. Furthermore, let $[i_1,j_1],[i_2,j_2],\ldots,[i_{n-k},j_{n-k}]$ be the matched pairs where 
$i_l<j_l$ for all $l$. Then we have:
\begin{enum} 
\item $c(i,j)$ is even if and only if $\pi_i+j$ is even.
\item $\inv(\pi)=c(i_1,j_1)+c(i_2,j_2)+\ldots+c(i_{n-k},j_{n-k})+n-k$. 
\end{enum}
\end{lem}  
\vspace*{-2ex}

\begin{proof}
{\bf a)} Since $\pi\in T_n$, we have $c_1(i,j)=0$. Moreover, there is no dot northeast of 
$(i,\pi_i)$ or southwest of $(j,\pi_j)$. Hence 
$c_2(i,j)+c_3(i,j)+2c_4(i,j)=n-\pi_i+n-j$. Thus we have $c(i,j)=c_2(i,j)+c_3(i,j)\equiv 
\pi_i+j\bmod2$.\\
{\bf b)} Since $\pi$ is bi-increasing, any inversion of $\pi$ has to 
be a pair $(i,j)$ with an excedance $i$ and an anti-excedance $j$.\\ 
First assume that $i$ is a matched excedance, and let $j$ be its match. 
Furthermore, let $(i,j')$ be an inversion. If $j\le j'$ we have 
$\pi_j\le\pi_{j'}<\pi_i$. There are exactly $c_3(i,j)+1$ such integers $j'$. 
In case $j'<j$, the anti-excedance $j'$ is 
matched with an excedance $i'$ for which $i'<i$. (Otherwise, $j'$ would be 
matched with $i$ by the definition of the matching.) Thus $(i,\pi_i)$ belongs 
to the dots which are counted by $c_2(i',j')$.\\        
Now let $i$ be an unmatched excedance. If $(i,j)$ is an inversion, then $j$ must be 
matched. (Otherwise, we could match $i$ with $j$.) By the construction, the 
match $i'$ of $j$ satisfies $i'<i$ (and hence $\pi_{i'}<\pi_i$). Thus $(i,\pi_i)$ is contained in 
the region whose dots are counted by $c_2(i',j)$. 
\end{proof}
By the definition, we have 
$\srs(\pi)=\pi_{i_1}+\pi_{i_2}+\ldots+\pi_{i_{n-k}}+j_1+j_2+\ldots+j_{n-k}$. 
Therefore, from the lemma we immediately obtain the assertion of Proposition \ref{sign}: 
a $321$-avoiding permutation is even if and only if the sum of the elements of 
the second row of the two tableaux, increased by the length of this row, is 
even.     
\vspace*{1cm}


\setcounter{section}{3}\setcounter{theo}{0}

\centerline{\large{\bf 3}\hspace*{0.25cm}
{\sc Proof of Theorem 1.1}}
\vspace*{0.5cm}

In this section, we construct a bijection on $T_n$ which proves the main 
result. Its essential part is an involution on the ballot sequences of length 
$n$ having a given number $k$ of $1$'s.\\[2ex] 
We define the {\it sign} of a ballot sequence $b$ to be $1$ if the sum of all integers $i$ with $b_i=-1$ is even and 
$-1$ otherwise, and write $\sg(b)$ to denote it.\\[2ex]
For a ballot sequence $b$, let $\varepsilon(b)$ be the smallest {\bf even} integer 
$i$ for which $b_i=-b_{i+1}$. If there is no such integer set 
$\varepsilon(b)=0$. By $A_{n,k}$ we denote the set of all ballot sequences $b$ of length $n$ 
having $k$ 1's and satisfying $\varepsilon(b)>0$, and by $A^*_{n,k}$ the set of all sequences $b$ 
with $\varepsilon(b)=0$.\\[2ex]
Consider now the following map $\phi$ on $A_{n,k}$. Let $b\in A_{n,k}$ satisfy $\varepsilon(b)=j$. Then define $c=\phi(b)$ to be the sequence 
which is obtained from $b$ by exchanging the elements $b_j$ and $b_{j+1}$, i.e., $c_j=-b_j$, 
$c_{j+1}=-b_{j+1}$, and $c_i=b_i$ otherwise.    

\begin{prop} \label{phi}
\begin{enum}
\item[]
\item The map $\phi$ is a sign-reversing involution on $A_{n,k}$. 
\item For odd $n$, we have 
$|A^*_{n,k}|=b(\frac{n-1}{2},\frac{k-1}{2})$ if $k$ is odd, and 
$|A^*_{n,k}|=0$ otherwise.
\item For even $n$, we have 
$|A^*_{n,k}|=b(\frac{n}{2}-1,\lfloor\frac{k-1}{2}\rfloor)$ for all $k$.
\end{enum}
\end{prop}     
\vspace*{-0.1cm}

\begin{proof}
{\bf a)} If $b\in A_{n,k}$ with $\varepsilon(b)=j$, then the sequence $c=\phi(b)$ 
belongs to $A_{n,k}$ as well. Note that $b_1+\ldots+b_{j-1}\ge 1$ since $j$ is 
even, and hence $c_1+\ldots+c_j=b_1+\ldots+b_{j-1}-b_j\ge0$. Obviously, $\phi$ is an involution. 
By the construction, the sum of the positions of the $-1$'s in $b$ and $c$ differs by 1. Hence 
$\sg(b)=-\sg(c)$.\\
{\bf b)}, {\bf c)} A sequence $b$ belongs to $A^*_{n,k}$ if and only if 
$b_{2i}=b_{2i+1}$ for $1\le i\le\lfloor\frac{n-1}{2}\rfloor$.\\ 
It is easy to see that there is a bijection between the sequences in $A^*_{n,k}$ and 
all the ballot sequences of length $\lfloor\frac{n-1}{2}\rfloor$ having a certain number of 
$1$'s. In case of odd $n$, set $d=b_2b_4\cdots b_{n-1}$. (Note that $b_2=b_3=1$; 
otherwise, we would have $b_1+b_2+b_3=-1$.) The number $k$ of 1's in $b$ is always 
an odd integer since $b_1=1$. Therefore, $d$ is a ballot sequence with $\frac{k-1}{2}$ 
components equal $1$. If $n$ is even, then the sequence $d$ is defined as $d=b_2b_4\cdots b_{n-2}$. 
Because we have $b_n=1$ if $k$ is even and $b_n=-1$ otherwise we may omit the final 
element. Hence the number of 1's in $d$ is equal to $\lfloor\frac{k-1}{2}\rfloor$.
\end{proof}
\vspace*{2ex}

It should be clear now how to construct the involution for proving Theorem \ref{main theorem}.\\ 
Given $\pi\in T_n$ with $\lis(\pi)=k$, let $p$ and $q$ be the ballot sequences defined by 
the tableaux $P$ and $Q$ resulting from the RSK algorithm. Define $\sigma=\Phi(\pi)$ to be the 321-avoiding permutation 
whose associated pair $(P',Q')$ of tableaux is given by the pair 
\bdpm
(p',q')=\left\{\ba{ccl}
(\phi(p),q)&&\mbox{if }p\in A_{n,k}\\
(p,\phi(q))&&\mbox{if }p\in A^*_{n,k}\mbox{ and }q\in A_{n,k}\\
(p,q)&&\mbox{otherwise}
\ea\right.
\edpm
of ballot sequences.\\[2ex]
Since $\phi$ preserves the number of 1's in a ballot sequence we have 
$\lis(\sigma)=\lis(\pi)$. By the definition, $\srs(\pi)$ is even if and only 
if $p$ and $q$ have the same sign. Consequently, $\Phi$ reverses the sign of 
the permutation if $p\in A_{n,k}$ or $q\in A_{n,k}$. (This follows from 
\ref{sign} and \ref{phi}a.)\\
Let now $\pi$ be a fixed point of $\Phi$, that is, $\varepsilon(p)=\varepsilon(q)=0$. 
Clearly, the number of such permutations $\pi$ is just $|A^*_{n,k}|^2$. 
In case of odd $n$ (and hence odd $k$), we have $\sg(\pi)=(-1)^{\srs(\pi)}$. 
By the proof of Proposition \ref{phi}b, $p_i=-1$ for any even integer $i$ if and only if 
$p_{i+1}=-1$. The same relation holds for the sequence $q$, too. Therefore, 
$\sum_{p_i=-1} i$ and $\sum_{q_i=-1} i$, respectively, are the sum of 
$\frac{n-k}{2}$ odd numbers. In particular, $p$ and $q$ are of the same sign. 
Thus $\srs(\pi)$ is even, and $\sg(\pi)=1$. By similar reasoning, 
we obtain $\sg(\pi)=(-1)^k$ if $n$ is even.\\[2ex]
Summarized, $\Phi$ yields the relations
\beas
&&e(2n+1,2k+1)-o(2n+1,2k+1)\;=\;b(n,k)^2\\
&&e(2n+1,2k)-o(2n+1,2k)\;=\;0\\
&&e(2n+2,2k+1)-o(2n+2,2k+1)\;=\;-b(n,k)^2\\
&&e(2n+2,2k+2)-o(2n+2,2k+2)\;=\;b(n,k)^2
\eeas
for all $n$ and $k$. Here let $e(n,k)$ and $o(n,k)$ denote the number of even and 
odd permutations in $T_n$, respectively, whose longest increasing subsequence is exactly 
of length $k$. 
\vspace*{1cm}


\setcounter{section}{4}\setcounter{theo}{0}

\centerline{\large{\bf 4}\hspace*{0.25cm}
{\sc A combinatorial proof of Adin-Roichman's result}}
\vspace*{0.5cm}

In \cite{adin-roichman}, Adin and Roichman proved the sign-balance on the set of $321$-avoiding 
permutations having a given last descent by a recursion formula for the 
generating function and asked for a combinatorial proof. By modifying $\Phi$ 
slightly, we obtain the desired involution. Moreover, this involution even proves 
a further refinement of our main result and Theorem \ref{adin}.\\[2ex]
For $\pi\in{\cal S}_n$, a {\it descent} of $\pi$ is an integer $i$ for 
which $\pi_i>\pi_{i+1}$. The set of the descents of $\pi$ we denote by ${\sf 
D}(\pi)$, its maximum element by $\ldes(\pi)$. (Define $\ldes({\rm id})=0$.)\\ 
Because of the twice order of $321$-avoiding permutations, 
an integer $i$ is a descent of $\pi\in T_n$ if and only if $i$ is an excedance but $i+1$ is none. In particular, the 
maximum descent is the maximum excedance as well.\\ 
By \cite[Cor. 3.13]{reifegerste}, there are 
$b(n+d-1,n-1)=\frac{n-d}{n+d}{n+d\choose d}$ permutations $\pi\in T_n$ 
for which $\ldes(\pi)=d$. 

\begin{theo}[{\cite[Thm. 4.1]{adin-roichman}}] \label{adin}
For all $n\ge1$, we have
\bdpm
\sum_{\pi\in T_{2n+1}} \sg(\pi)\cdot q^{\ldes(\pi)}\:= 
\sum_{\pi\in T_n} q^{2\ldes(\pi)}\;,\quad
\sum_{\pi\in T_{2n}} \sg(\pi)\cdot q^{\ldes(\pi)}\:=
(1-q)\sum_{\pi\in T_n} q^{2\ldes(\pi)}.
\edpm
\end{theo}

By the description of the RSK algorithm, the last descent is just the maximum element $i$ in the 
first row of the tableau $Q$ for which $i+1$ appears in the second row of $Q$. 
(Moreover, an integer $i$ is a descent of $\pi\in T_n$ if and only if $i$ 
appears in the first row while $i+1$ is in the second row of $Q$.) Thus the bijection $\Phi$ also preserves 
the last descent of $\pi$ if $\varepsilon(p)>0$ or 
$\varepsilon(q)\le\ldes(\pi)-2$. Consequently, the involution we look for will nearly coincide with $\Phi$; we only have to modify the cases
$\varepsilon(p)=0,\;\varepsilon(q)\ge\ldes(\pi)-1$ and 
$\varepsilon(p)=\varepsilon(q)=0$.\\[2ex]
For a ballot sequence $b$, let $\delta(b)$ be the greatest integer $i$ for which 
$b_i=-b_{i+1}=1$. (As mentioned above, we have $\delta(q)=\ldes(\pi)$ if $q$ is the 
sequence defined by the insertion tableau of $\pi$.) 
In addition to $A_{n,k}$ and $A^*_{n,k}$, we use $B_{n,k}$, $B^*_{n,k}$, and 
$B^\times_{n,k}$ to denote the set of all ballot sequences $b$ of length $n$ 
having $k$ 1's and satisfying $0<\varepsilon(b)<\delta(b)-1$, 
$\varepsilon(b)=\delta(b)-1$, and $\varepsilon(b)\ge\delta(b)$, respectively. 
(The set $A_{n,k}$ is just the union of $B_{n,k}$, $B^*_{n,k}$, and $B^\times_{n,k}$.)

\begin{lem} \label{psi}
Let $\pi\in T_n$ be a permutation for which $\lis(\pi)=k$ and $\ldes(\pi)=d$, and $(P,Q)$ its associated pair of tableaux. Furthermore, let 
$p$ and $q$ be the ballot sequences defined by $P$ and $Q$, respectively. Assume that 
$p\in A_{n,k}^*$ and $q\in A^*_{n,k}\cup B^*_{n,k}\cup B^\times_{n,k}$. Then we have:
\begin{enum} 
\item If $d$ is even, then $\pi$ is even.
\item If both $n$ and $d$ are odd, then $\pi$ is even if and only if $q\in 
A_{n,k}^*$. Moreover, there are as many sequences $q\in A_{n,k}^*$ as sequences $q\in 
B_{n,k}^*$ satisfying $\delta(q)=d$.
\item If $n$ is even and $d$ is odd, then $\pi$ is even if and only if $q\in 
A_{n,k}^*$ and $k$ is even. Moreover, for any even $k$ there are as many sequences $q\in A_{n,k}^*$ as sequences $q\in 
B_{n,k}^*$ with $q_n=1$ satisfying $\delta(q)=d$.
\end{enum}
\end{lem} 
\vspace*{-2ex}
 
\begin{proof} 
By the definition of $\delta$, we have $q=q_1q_2\cdots q_{d-1}1(-1)(-1)^a1^b$ where 
$a,b\ge0$. (Here exponentiation denotes repetition.)\\
{\bf a)} If $d$ is even, then we have $\varepsilon(q)=d$. In particular, all 
the components $q_i=-1$ with $i<d$ appear in pairs. 
First let $n$ be odd. By Proposition \ref{phi}b, then $k$ must be odd as well 
since $p\in A^*_{n,k}$. As already discussed in the previous section, $\sum_{p_i=-1} i$ is the sum of 
$\frac{n-k}{2}$ odd numbers. Since $q$ contains the element $-1$ exactly $n-k$ 
times, the exponent $a$ is odd. Thus $\sum_{q_i=-1} i$ is also the sum of $\frac{n-k}{2}$ odd numbers. 
Hence $p$ and $q$ have the same sign, and $\srs(\pi)$ is even. In case of 
even $n$, we have to distinguish whether $k$ is even or not. By reasoning in a similar 
way as done for odd $n$, we can show that $\srs(\pi)$ is even if and only if $k$ is 
even. Thus $\sg(\pi)=1$.\\   
{\bf b)} Let both $n$ (and hence $k$) and $d$ be odd. By the definition, we have 
$q\in A^*_{n,k}$ if $q_{d-1}=1$, and $q\in B^*_{n,k}$ if $q_{d-1}=-1$. (Note 
that the exponent $a$ must be odd in the first case.)\\
Consider the following map $\psi$ that takes $q\in A^*_{n,k}$ to a sequence 
$q'\in B^*_{n,k}$. For $q\in A^*_{n,k}$, let $j$ be the last 
position of $-1$ in $q$. Define $q'$ to be the sequence which is obtained from $q$ by exchanging the elements $q_{d-1}$ and 
$q_j$. Since $d-1$ is even, we have $q_1+\ldots+q_{d-2}\ge 1$, and the exchange 
is permitted. Obviously, $q'\in B^*_{n,k}$ and $\delta(q')=d$. (Note that $j$ 
must be odd. Thus we may assume that $j>d+1$.) Moreover, $\psi$ reverses the 
sign. It is easy to see that $\psi$ is bijective. The inverse map just takes 
$q\in B^*_{n,k}$ to the sequence $q'\in A^*_{n,k}$ which is defined by 
$q'_{d-1}=q_j=1,\;q'_j=q_{d-1}=-1$ and $q'_i=q_i$ otherwise where $j$ is the 
position of the first 1 to the right of $q_d$.\\ 
Similarly to part a), we can prove that the signs of $p$ and $q$ coincide if $q\in A^*_{n,k}$. 
Applying $\psi$ yields the conversion.\\
{\bf c)} Let $n$ be even now and $d$ odd. For even $k$ we can prove the 
assertion by the same arguments as used in part b). Note that the map 
$\psi$ is not defined if $q\in B^*_{n,k}$ ends with an element $-1$. 
In particular, we have $\sg(p)=\sg(q)$ (and hence $\sg(\pi)=1$) if and only if 
$q\in A^*_{n,k}$. If $k$ is odd, then $p$ and $q$ are of the same sign. On the 
one side, we have $p_n=-1$, on the other side either there is an even integer $i$ for 
which $q_i=-1$ and $q_{i+1}=1$ or $d=n-1$. Hence $\sg(\pi)=-1$. 
\end{proof}
\vspace*{1ex}

Let $\pi\in T_n$ satisfy $\lis(\pi)=k$ and $\ldes(\pi)=d$. 
As before, let $p$ and $q$ be the ballot sequences defined by the tableaux $P$ 
and $Q$ obtained from the RSK correspondence. Define $\sigma=\Psi(\pi)$ to be the 321-avoiding permutation 
whose associated tableaux $P'$ and $Q'$ are given by the sequences
\bdpm
(p',q')=\left\{\ba{ccl}
(\phi(p),q)&&\mbox{if }p\in A_{n,k}\\
(p,\phi(q))&&\mbox{if }p\in A^*_{n,k}\mbox{ and }q\in B_{n,k}\\
(p,\psi(q))&&\mbox{if }p\in A^*_{n,k}\mbox{ and }q\in A^*_{n,k}\cup B^{**}_{n,k};\; 
n-k\mbox{ is even and }d\mbox{ is odd}\\
(p,q)&&\mbox{otherwise}
\ea\right.
\edpm
where $\psi$ is defined in the proof of Lemma \ref{psi}b, and $B^{**}_{n,k}$ 
denotes the set of all sequences $q\in B^*_{n,k}$ for which $q_n=1$. By the previous 
discussion, $\Psi$ is an involution on $T_n$ which preserves both $\lis$ and $\ldes$, and reverses the sign if $\pi$ is not a fixed point.\\[2ex] 
Finally, we enumerate the fixed points of $\Psi$ regarding the last descent.  

\begin{prop}
Let $f(n,d)$ be the number of permutations $\pi\in T_n$ for which 
$\Psi(\pi)=\pi$ and $\ldes(\pi)=d$.
\begin{enum} 
\item For odd $n$, we have $f(n,d)=b(\frac{n+d-3}{2},\frac{n-3}{2})$ if $d$ is even, and $f(n,d)=0$ otherwise.
\item For even $n$, we have $f(n,d)=b(\lfloor\frac{n+d-2}{2}\rfloor,\frac{n-2}{2})$ for all $d$.
\end{enum}
\end{prop}
\vspace*{-2ex} 

\begin{proof}
{\bf a)} First let $n$ be odd. In this case, $\pi$ is a fixed point if and 
only if $p\in A^*_{n,k}$ and $\varepsilon(q)=d$. (On condition that $d$ is 
even, this is the only possibility for $q\in A^*_{n,k}\cup B^*_{n,k}\cup 
B^\times_{n,k}$.) As shown in the proof of Proposition \ref{phi}b, the sequence 
$p$ corresponds in a one-to-one fashion to a ballot sequence of length 
$\frac{n-1}{2}$ with $\frac{k-1}{2}$ components equal 1. (Note that $k$ has to 
be odd.) On the other side, there is a bijection between ballot sequences 
$q$ of length $n$ with an odd number $k$ of 1's satisfying 
$\varepsilon(q)=\delta(q)=d$ and ballot sequences $q'$ of length $\frac{n-1}{2}$ 
having $\frac{k-1}{2}$ 1's and satisfying $\delta(q')=\frac{d}{2}$. Set $q'=q_2q_4\cdots 
q_{d-2}q_dq_{d+2}\cdots q_{n-1}$. Note that we have $q_i=q_{i+1}$ for all even $i<d$. Furthermore, 
$q_i=-1$ for $d+1\le i\le j$ and $q_i=1$ for $j+1\le i\le n$ where $j$ must be an 
even integer because $n-k$ is even. In particular, $q_{d+2}=-1$. Hence we have 
$\delta(q')=\frac{d}{2}$. Consequently, $(p,q)$ corresponds to a pair of ballot 
sequences that is associated with a permutation $\tau\in T_{\frac{n-1}{2}}$ for 
which $\ldes(\tau)=\frac{d}{2}$ (and $\lis(\tau)=\frac{k-1}{2}$). Thus $f(n,d)=b(\frac{n-1}{2}+\frac{d}{2}-1,\frac{n-1}{2}-1)$.\\
{\bf b)} Now consider the case of even $n$. If $d$ is even, we have $\varepsilon(q)=d$ 
again. For even $k$, the situation is similar to part a). The sequence $q$ 
corresponds to a ballot sequence $q'$ of length $\frac{n}{2}$ with 
$\frac{k}{2}$ elements equal 1 and $\delta(q')=\frac{d}{2}$; set 
$q'=q_2q_4\cdots q_dq_{d+2}\cdots q_n$. If $k$ is odd, we can not assume that 
$q_{d+2}=-1$. Thus we must add the element $q_{d+1}$ to $q'$ to obtain 
$\delta(q')=\frac{d}{2}$. On the other 
side, we may omit the final element since we have $q_n=1$. Hence $q$ can be 
identified with $q'=q_2q_4\cdots q_dq_{d+1}q_{d+2}q_{d+4}\cdots q_{n-2}$; a ballot sequence of length $\frac{n}{2}$ with 
$\frac{k-1}{2}$ 1's and $\delta(q')=\frac{d}{2}$. 
For any fixed point $\pi$ of $\Psi$, the sequence $p$ belongs to $A^*_{n,k}$. By Proposition 
\ref{phi}c, thus $p$ corresponds to a ballot 
sequence of length $\frac{n}{2}-1$ with $\lfloor\frac{k-1}{2}\rfloor$ 1's. 
Define $p'$ to be the sequence obtained from $p$ by this bijection and adding 
an element 1 if $k$ is even and $-1$ otherwise. This yields a one-to-one correspondence between 
$(p,q)$ and a pair $(p',q')$ of ballot sequences of length $\frac{n}{2}$ with 
$\lfloor\frac{k}{2}\rfloor$ 1's satisfying $\delta(q')=\frac{d}{2}$. 
Consequently, $f(n,d)$ counts the number of permutations in $T_{\frac{n}{2}}$ 
having the last descent $\frac{d}{2}$.\\ 
It is not difficult to see that for any even $d$ there are  
as many sequences $q$ with $\delta(q)=d$ as such ones with $\delta(q)=d+1$ among the 
ballot sequences arising from the insertion tableau of a fixed point of $\Psi$. 
\end{proof}

By Lemma \ref{psi}, the fixed points of $\Psi$ are odd permutations if and only 
if $n$ is even and $d$ is odd. This completes the proof of Theorem \ref{adin}.\\[2ex]   
As a by-product we obtain the following refinement of the theorems \ref{main theorem} and 
\ref{adin}.

\begin{cor} 
For all $n\ge1$, we have
\beas
\sum_{\pi\in T_n} q^{2\lis(\pi)+1}t^{2\ldes(\pi)}&=&
\sum_{\pi\in T^*_{2n+1}} \sg(\pi)\cdot q^{\lis(\pi)}t^{\ldes(\pi)}\\[1ex]
\sum_{\pi\in T_n} q^{2\lis(\pi)+1}t^{2\ldes(\pi)}&=&
\sum_{\pi\in T^*_{2n}} \sg(\pi)\cdot q^{\lis(\pi)}t^{\ldes(\pi)}+
q\sum_{\pi\in T^\star_{2n}} \sg(\pi)\cdot q^{\lis(\pi)}t^{\ldes(\pi)}.
\eeas
where $T^*_n$ (resp. $T^\star_n$) denotes the set of all permutations $\pi\in T_n$ 
whose last descent is even and whose longest increasing subsequence is of odd 
(resp. even) length. 
\end{cor}
\vspace*{0.5cm}


\setcounter{section}{5}\setcounter{theo}{0}

\centerline{\large{\bf 5}\hspace*{0.25cm}
{\sc Final remarks}}
\vspace*{0.5cm}

In their paper, Adin and Roichman investigate a further permutation statistic 
on $T_n$: the position of the letter $n$ in $\pi$, denoted by $\lind$. Using generating functions 
they prove the equidistribution of $\ldes$ and $\lind-1$. We give a simple 
bijective proof.  

\begin{theo}[{\cite[Thm. 3.1]{adin-roichman}}]
The statistics $\ldes$ and $\lind-1$ are equidistributed over $T_n$. Moreover, for any 
$B\subseteq[n-2]$ they are equidistributed 
over the set
\bdpm
\{\pi\in T_n:{\sf D}(\pi^{-1})\cap[n-2]=B\}
\edpm
where $[n]$ means the set $\{1,\ldots,n\}$.
\end{theo}

\begin{proof}
Let $\pi\in T_n$ be a permutation for which $\ldes(\pi)=d$. Define 
$\sigma$ to be the permutation which is obtained from $\pi$ by deleting the 
letter $n$, followed by inserting the element $n$ between the positions $d$ and 
$d+1$.\\
Obviously, this map is a bijection on $T_n$ since the construction 
of $\sigma$ preserves the order of the excedance letters and non-excedance letters, 
respectively. Clearly, $\lind(\sigma)=d+1$.\\
Furthermore, we have ${\sf D}(\pi^{-1})\cap[n-2]={\sf 
D}(\sigma^{-1})\cap[n-2]$. By the definition, the descent set ${\sf D}(\pi^{-1})$ of the inverse of 
$\pi$ contains all the integers $i$ for which the element $i+1$ appears to the left of 
the element $i$ in $\pi$. (If $n-1\in{\sf D}(\pi^{-1})$, then we have 
$\ldes(\pi)=\lind(\pi)=d$. In case $\pi_{d+1}=n-1$, the set ${\sf D}(\sigma^{-1})$ 
does not contain $n-1$.) 
\end{proof}

Consequently, we also obtain sign-balance on $T_n$ concerning the statistic $\lind$.
\vspace*{1.5cm}

{\bf Acknowledgement} I would like to thank Herb Wilf for the interest taken in this work.
\newpage


\centerline{\large\sc References}
\vspace*{0.5cm}

\begin{enumbib}

\bibitem{adin-roichman}
R. M. Adin and Y. Roichman,
Equidistribution and Sign-Balance on $321$-Avoiding Permutations, 
preprint, 2003, math.CO/0304429.

\bibitem{deutsch etal}
E. Deutsch, A. J. Hildebrand, and H. S. Wilf,
Longest increasing subsequences in pattern-restricted permutations, 
preprint, 2003, math.CO/0304126.

\bibitem{elizalde-pak}
S. Elizalde and I. Pak,
Bijections for refined restricted permutations, 
preprint, 2002, math.CO/0212328.

\bibitem{feller}
W. Feller, 
{\it An Introduction to Probability Theory and Its Applications}, vol. I,
Wiley, New York, 1968.

\bibitem{reifegerste}
A. Reifegerste,
The excedances and descents of bi-increasing permutations, 
preprint, 2002, math.CO/\\0212247.

\bibitem{simion-schmidt}
R. Simion and F. W. Schmidt, 
Restricted Permutations, 
{\it Europ. J. Combinatorics} {\bf 6} (1985), 383-406. 

\bibitem{stanley}
R. P. Stanley, 
{\it Enumerative Combinatorics}, vol. II,
Cambridge University Press, 1999.

\bibitem{stanley1}
R. P. Stanley, 
Some Remarks on Sign-Balanced and Maj-Balanced Posets, 
preprint, 2002, math.CO/\\0211113.

\end{enumbib}

\end{document}